\documentclass[12pt]{article}
\usepackage{latexsym, amssymb, amsmath, amscd, amsfonts, epsfig, graphicx, colordvi,verbatim,ifpdf,extarrows}
\usepackage{amsfonts, amsmath, amssymb}
\usepackage{amssymb,amsfonts,amsmath,latexsym,epsfig,cite, psfrag,eepic,color}
\usepackage{amscd,graphics}
\usepackage{latexsym, amssymb,  amsmath,amscd, amsfonts, epsfig, graphicx, colordvi,amsthm}

\usepackage{graphicx}
\usepackage{epstopdf}
\usepackage{color}
\usepackage{ifpdf}
\usepackage{fancybox}
\usepackage[font=small,labelfont=bf,labelsep=none]{caption}
\usepackage{float}

\allowdisplaybreaks

\newtheorem{thm}{Theorem}[section]

\newtheorem{defi}[thm]{Definition}
\newtheorem{lem}[thm]{Lemma}
\newtheorem{core}[thm]{Corollary}

\def\pf{\noindent{\it Proof.} }
\def\qed{\nopagebreak\hfill{\rule{4pt}{7pt}}
\medbreak}

\numberwithin{equation}{section}

\def\qed{\nopagebreak\hfill{\rule{4pt}{7pt}}
\medbreak}

\newlength{\boxedparwidth}
  {\begin{center} \begin{tabular}{|@{\hspace{.315in}}c@{\hspace{.15in}}|}
                  \hline \\ \begin{minipage}[t]{\boxedparwidth}
                  \setlength{\parindent}{.25in}}%
  {\end{minipage} \\ \\ \hline \end{tabular} \end{center}}

\begin{document}

\begin{center}
{\Large \bf Overpartitions with separated overlined parts and non-overlined parts}
\end{center}

\begin{center}
 {Y.H. Chen}$^{1}$, {Thomas Y. He}$^{2}$, {Y. Hu}$^{3}$ and
  {Y.X. Xie}$^{4}$ \vskip 2mm

$^{1,2,3,4}$ School of Mathematical Sciences, Sichuan Normal University, Chengdu 610066, P.R. China

   \vskip 2mm

  $^1$chenyh@stu.sicnu.edu.cn, $^2$heyao@sicnu.edu.cn,  $^3$huyue@stu.sicnu.edu.cn,  $^4$xieyx@stu.sicnu.edu.cn
\end{center}

\vskip 6mm   {\noindent \bf Abstract.} Recently, Andrews considered the partitions with parts separated by parity, in which parts of a given parity are all smaller than those of the other parity. Inspired from the partitions with parts separated by parity, we investigate the overpartitions with separated overlined parts and non-overlined parts, in which the sizes of overlined parts (resp. non-overlined parts) are greater than or equal to those of non-overlined parts (resp. overlined parts).

\noindent {\bf Keywords}: overpartitions, separable overpartition classes, distinct partitions, overlined parts, non-overlined parts

\noindent {\bf AMS Classifications}: 05A17, 11P83

\section{Introduction}

A partition $\pi$ of a positive integer $n$ is a finite non-increasing sequence of positive integers $\pi=(\pi_1,\pi_2,\ldots,\pi_m)$ such that $\pi_1+\pi_2+\cdots+\pi_m=n$. The empty sequence forms the only partition of zero. The $\pi_i$ are called the parts of $\pi$.  Let $\ell(\pi)$ be the number of parts of $\pi$. The weight of $\pi$ is the sum of parts, denoted $|\pi|$.  We use $\mathcal{P}(n)$ to denote the set of all partitions of $n$.
The conjugate $\pi'$ of a partition $\pi$ is a partition such that the $i$-th part $\pi'_i$ is the number of parts of $\pi$ greater than or equal to $i$.

In \cite{Andrews-2018,Andrews-2019}, Andrews considered the partitions with parts separated by parity, in which parts of a given parity are all smaller than those of the other parity. In \cite{Andrews-2022}, Andrews introduced separable integer partition classes and analyzed some well-known theorems. Passary \cite[Section 3.2]{Passary-2019} and Chen, He, Tang and Wei  \cite[Section 3]{Chen-He-Tang-Wei-2024} studied the partitions with parts separated by parity from the view of the point of
 separable integer partition classes.

Inspired from the partitions with parts separated by parity, we consider the overpartitions with separated overlined parts and non-overlined parts. An overpartition, introduced by Corteel and Lovejoy \cite{Corteel-Lovejoy-2004},  is a partition such that the first occurrence of a number can be overlined.
For example, there are fourteen overpartitions of $4$:
\[(4),(\overline{4}),(3,1),(\overline{3},1),(3,\overline{1}),(\overline{3},\overline{1}),(2,2),(\overline{2},2),
\]
\[(2,1,1),(\overline{2},1,1),(2,\overline{1},1),(\overline{2},\overline{1},1),(1,1,1,1),(\overline{1},1,1,1).\]

We will use $\overline{\mathcal{P}}(n)$ to denote the set of all overpartitions of $n$. Let $\pi$ be an overpartition. We sometimes write  $\pi=\left(1^{f_1(\pi)}\bar{1}^{f_{\bar{1}}(\pi)}2^{f_2(\pi)}\bar{2}^{f_{\bar{2}}(\pi)}\cdots\right)$, where $f_t(\pi)$ (resp. $f_{\bar{t}}(\pi)$) denotes the number of parts equal to $t$ (resp. $\bar{t}$) in $\pi$.
For a part $\pi_i$ of $\pi$, we say that $\pi_i$ is of size $t$ if $\pi_i=t$ or $\overline{t}$, denoted $|\pi_i|=t$. We define $t+ b$ (resp. $\overline{t}+ b$) as a non-overlined part (resp. an overlined part) of size  $t+b$.

In this article, we investigate four types of overpartitions.
\begin{itemize}
\item Let $\mathcal{G}_N^{O}$ (resp. $\mathcal{E}_N^{O}$) be the set of overpartitions such that the sizes of overlined parts (if exist) are greater than (resp. greater than or equal to) those of non-overlined parts.

\item Let $\mathcal{G}_O^{N}$ (resp. $\mathcal{E}_O^{N}$) be the set of overpartitions such that non-overlined parts must appear and the sizes of non-overlined parts are greater than (resp. greater than or equal to) those of overlined parts.
\end{itemize}

For example, if we consider the overpartitions of $4$, then we have
\[(4),(\overline{4}),(3,1),(\overline{3},1),(\overline{3},\overline{1}),(2,2),(2,1,1),(\overline{2},1,1),(1,1,1,1)\in\mathcal{G}_N^{O};\]
\[(4),(\overline{4}),(3,1),(\overline{3},1),(\overline{3},\overline{1}),(2,2),(\overline{2},2),(2,1,1),(\overline{2},1,1),(\overline{2},\overline{1},1),(1,1,1,1),(\overline{1},1,1,1)\in\mathcal{E}_N^{O};\]
\[(4),(3,1),(3,\overline{1}),(2,2),(2,1,1),(1,1,1,1)\in\mathcal{G}_O^{N};\]
\[(4),(3,1),(3,\overline{1}),(2,2),(\overline{2},2),(2,1,1),(2,\overline{1},1),(1,1,1,1),(\overline{1},1,1,1)\in\mathcal{E}_O^{N}.\]

For easier expression, we use the following notations related to an overpartition $\pi$.
\begin{itemize}
\item Let $LN(\pi)$ (resp. $SN(\pi)$) be the size of the largest (resp. smallest) non-overlined part of $\pi$ if there exist non-overlined parts in $\pi$, and $LN(\pi)=0$ (resp. $SN(\pi)=0$) otherwise.

\item Let $LO(\pi)$ (resp. $SO(\pi)$) be the size of the largest (resp. smallest) overlined part of $\pi$ if there exist overlined parts in $\pi$, and $LO(\pi)=0$ (resp. $SO(\pi)=0$) otherwise.

\item Let $\ell_{N\geq O}(\pi)$ (resp. $\ell_{N>O}(\pi)$) be the number of non-overlined parts of size greater than or equal to (resp. greater than) $SO(\pi)$ in $\pi$.

\item Let $\ell_{O\geq N}(\pi)$ (resp. $\ell_{O>N}(\pi)$) be the number of overlined parts of size greater than or equal to (resp. greater than) $SN(\pi)$ in $\pi$.

\end{itemize}

Clearly, $\mathcal{G}_N^{O}$ (resp. $\mathcal{E}_N^{O}$) is the set of overpartitions $\pi$ such that if $SO(\pi)\geq 1$ then $SO(\pi)>LN(\pi)$ (resp. $SO(\pi)\geq LN(\pi)$), and  $\mathcal{G}_O^{N}$ (resp. $\mathcal{E}_O^{N}$) is the set of overpartitions $\pi$ such that $SN(\pi)\geq 1$ and $SN(\pi)>LO(\pi)$ (resp. $SN(\pi)\geq LO(\pi)$).

For a nonempty overpartition $\pi$, by definition, we see that  $\ell_{N\geq O}(\pi)=0$ if $\pi\in\mathcal{G}_N^{O}$ and $SO(\pi)\geq 1$; $\ell_{N>O}(\pi)=0$ if $\pi\in\mathcal{E}_N^{O}$ and $SO(\pi)\geq 1$; $\ell_{O\geq N}(\pi)=0$ if $\pi\in\mathcal{G}_O^{N}$; $\ell_{O>N}(\pi)=0$ if $\pi\in\mathcal{E}_O^{N}$.

We will show that $\mathcal{G}_N^{O}$, $\mathcal{E}_N^{O}$, $\mathcal{G}_O^{N}$ and $\mathcal{E}_O^{N}$ are separable overpartition classes introduced by Chen, He, Tang and Wei \cite[Definition 4.1]{Chen-He-Tang-Wei-2024}, which is an extension of separable integer partition classes given by Andrews \cite{Andrews-2022}.

\begin{defi}[separable overpartition classes]
A separable overpartition class $\mathcal{P}$ is a set of overpartitions satisfying the following{\rm:}

There is a subset $\mathcal{B}\subset\mathcal{P}$ {\rm(}$\mathcal{B}$ is called the basis of $\mathcal{P}${\rm)} such that for each integer $m\geq 1$, the number of overpartitions in $\mathcal{B}$ with $m$ parts is finite and every overpartition in $\mathcal{P}$ with $m$ parts is uniquely of the form
\begin{equation}\label{over-form-1}
(b_1+\pi_1)+(b_2+\pi_2)+\cdots+(b_m+\pi_m),
\end{equation}
where $(b_1,b_2,\ldots,b_m)$ is an overpartition in $\mathcal{B}$ and $(\pi_1,\pi_2,\ldots,\pi_m)$ is a  non-increasing sequence of nonnegative integers. Moreover, all overpartitions of the form \eqref{over-form-1} are in $\mathcal{P}$.
\end{defi}

For $m\geq 1$, let $b_\mathcal{B}(m)$ be the generating function for the overpartitions in $\mathcal{B}$ with  $m$ parts. The generating function for  the overpartitions in $\mathcal{P}$ with  $m$ parts is
 \[\frac{b_\mathcal{B}(m)}{(q;q)_m}.\]

 Here and in the sequel, we assume that $|q|<1$ and use the standard notation \cite{Andrews-1976}:
\[(a;q)_\infty=\prod_{i=0}^{\infty}(1-aq^i)\quad\text{and}\quad(a;q)_n=\frac{(a;q)_\infty}{(aq^n;q)_\infty}.\]

Motivated by the works of Kim, Kim and Lovejoy \cite{Kim-Kim-Lovejoy-2021} and Lin and Lin \cite{Lin-Lin-2024}, we consider the following partition functions.
 \begin{itemize}
\item[(1)] Let ${A}_{N\geq O}(n)$ (resp. ${B}_{N\geq O}(n)$) be the number of overpartitions $\pi$ of $n$ with $\ell_{N\geq O}(\pi)$ being even (resp. odd).

\item[(2)] Let ${A}_{N>O}(n)$ (resp. ${B}_{N>O}(n)$) be the number of overpartitions $\pi$ of $n$ with $\ell_{N>O}(\pi)$ being even (resp. odd).

\item[(3)] Let ${A}_{O\geq N}(n)$ (resp. ${B}_{O\geq N}(n)$) be the number of overpartitions $\pi$ of $n$ such that
  $SN(\pi)\geq 1$ and $\ell_{O\geq N}(\pi)$ is even (resp. odd).

\item[(4)] Let ${A}_{O>N}(n)$ (resp. ${B}_{O>N}(n)$) be the number of overpartitions $\pi$ of $n$ such that $SN(\pi)\geq 1$ and $\ell_{O>N}(\pi)$ is even (resp. odd).
\end{itemize}

This article is organized as follows. We will list the results of this article in Section 2.  In Section 3, we recall some necessary identities and give an involution on $\overline{\mathcal{P}}(n)$. Then, we will show the results in Sections 2.1 and 2.2 in Sections 4 and 5 respectively.

\section{The results of this article}

In this section, we will list the results related to $\mathcal{G}_N^{O}$ and $\mathcal{E}_N^{O}$ in Section 2.1 and the results related to $\mathcal{G}_O^{N}$ and $\mathcal{E}_O^{N}$ in Section 2.2.

\subsection{Results related to $\mathcal{G}_N^{O}$ and $\mathcal{E}_N^{O}$}
In view of separable overpartition classes, we can get
\begin{equation}\label{gen-G-O-N}
\sum_{\pi\in\mathcal{G}_N^{O}}q^{|\pi|}=(-q;q)_\infty^2,
\end{equation}
and
\begin{equation}\label{gen-E-O-N}
\sum_{\pi\in\mathcal{E}_N^{O}}q^{|\pi|}=(-q;q)_\infty^2+\frac{1}{(q;q)_\infty}-(-q;q)_\infty.
\end{equation}

We will give another two proofs of \eqref{gen-G-O-N} and \eqref{gen-E-O-N} in terms of the largest non-overlined part and the smallest overlined part. For $n\geq 0$, let ${G}_N^{O}(n)$ be the number of overpartitions of $n$ in  $\mathcal{G}_N^{O}$ and let $D_2(n)$ be the number of pairs of distinct partitions $(\alpha,\beta)$ such that $|\alpha|+|\beta|=n$. As a corollary of \eqref{gen-G-O-N}, we get
\begin{core}\label{core-g-d2}
For $n\geq0$,
\[{G}_N^{O}(n)=D_2(n).\]
\end{core}

In  \cite{Andrews-Newman-2019}, Andrews and Newman undertook a combinatorial study of the minimal excludant of a partition, which was earlier introduced by Grabner and Knopfmacher \cite{Grabner-Knopfmacher-2006} under the name ``smallest gap". The minimal excludant of a partition $\pi$ is the smallest positive integer that is not a part of $\pi$, denoted  $mex(\pi)$. For $n\geq 0$, Andrews and Newman defined
\[\sigma mex(n)=\sum_{\pi\in\mathcal{P}(n)}mex(\pi).\]
 They gave two proofs of the
following theorem.
\begin{thm}\cite{Andrews-Newman-2019}\label{simgamex-dn}
For $n\geq0$,
\[\sigma mex(n)=D_2(n).\]
\end{thm}
Ballantine and Merca \cite{Ballantine-Merca-2021} presented a combinatorial proof of Theorem \ref{simgamex-dn}.
By Corollary \ref{core-g-d2} and Theorem \ref{simgamex-dn}, we get the following corollary.
\begin{core}\label{core-g-sigma}
For $n\geq0$,
\[{G}_N^{O}(n)=\sigma mex(n).\]
\end{core}

For $n\geq 1$, let ${E}_{ON}(n)$ denote the number of overpartitions of $n$ such that the size of the smallest overlined part equals the size of the largest non-overlined part, that is, ${E}_{ON}(n)$ is the number of overpartitions $\pi$ of $n$ with $SO(\pi)=LN(\pi)$. For example, there are seven overpartitions counted by
${E}_{ON}(6)$, and so we have ${E}_{ON}(6)=7$.
\[(\overline{4},\overline{1},1),(\overline{3},3),(\overline{3},\overline{1},1,1),(\overline{2},2,2),(\overline{2},2,1,1),(\overline{2},\overline{1},1,1,1),(\overline{1},1,1,1,1,1).\]

It follows from \eqref{gen-G-O-N} and \eqref{gen-E-O-N} that
\begin{equation}\label{gen-E-O=N}
\sum_{n\geq 1}{E}_{ON}(n)q^{n}=\frac{1}{(q;q)_\infty}-(-q;q)_\infty.
\end{equation}

We will give another analytic proof of \eqref{gen-E-O=N}. For $n\geq 1$, let $R(n)$ denote the number of partitions of $n$ with repeated parts. For example, we have $R(6)=7$ since there are seven partitions of $6$ with repeated parts.
\[(4,1,1),(3,3),(3,1,1,1),(2,2,2),(2,2,1,1),(2,1,1,1,1),(1,1,1,1,1,1).\]

Clearly, the right hand side of \eqref{gen-E-O=N} is the generating function of $R(n)$. So, we have
\begin{core}\label{core-ON=R}
For $n\geq1$,
\[{E}_{ON}(n)=R(n).\]
\end{core}

Now, we turn to ${A}_{N\geq O}(n)$, ${B}_{N\geq O}(n)$, ${A}_{N>O}(n)$ and ${B}_{N>O}(n)$.
\begin{thm}\label{EG-OG-THM} For $n\geq 1$,
\begin{equation*}\label{EG-OG-eqn-1}
{A}_{N\geq O}(n)-{B}_{N\geq O}(n)\equiv0\pmod 2,
\end{equation*}
and
\begin{equation*}\label{EG-OG-eqn-2}
{A}_{N\geq O}(n)-{B}_{N\geq O}(n)\geq 0\text{ with strict inequality if }n\geq 3.
\end{equation*}

\end{thm}

Let $p_e(n)$ denote the number of partitions of $n$ with an even number of parts.
\begin{thm}\label{EE-OE-THM}
For $n\geq 1$,
\begin{equation*}\label{EE-OE-eqn-1}
{A}_{N>O}(n)-{B}_{N>O}(n)=2p_e(n).
\end{equation*}
\end{thm}
It yields that
\begin{core} For $n\geq 1$,
\begin{equation*}
{A}_{N>O}(n)-{B}_{N>O}(n)\equiv0\pmod 2,
\end{equation*}
and
\begin{equation*}
{A}_{N>O}(n)-{B}_{N>O}(n)\geq 0\text{ with strict inequality if }n\geq 2.
\end{equation*}
\end{core}

In Section 4, we will prove \eqref{gen-G-O-N}, \eqref{gen-E-O-N}, Theorems \ref{EG-OG-THM} and \ref{EE-OE-THM}, and we will present an analytic proof \eqref{gen-E-O=N} and bijective proofs of Corollary \ref{core-g-d2}, \ref{core-g-sigma} and \ref{core-ON=R}.

\subsection{Results related to $\mathcal{G}_O^{N}$ and $\mathcal{E}_O^{N}$}

By virtue of separable overpartition classes, we can get
\begin{equation}\label{gen-G-N-O}
\sum_{\pi\in\mathcal{G}_O^{N}}q^{|\pi|}=\sum_{m\geq1}\frac{1}{(q;q)_m} \sum_{k=0}^{m-1}q^{(k+1)m-\binom{k+1}{2}},
\end{equation}
and
\begin{equation}\label{gen-E-N-O}
\sum_{\pi\in\mathcal{E}_O^{N}}q^{|\pi|}=\sum_{m\geq1}\frac{1}{(q;q)_m} \left[q^m+\sum_{k=1}^{m-1}q^{km-\binom{k}{2}}\right].
\end{equation}

By considering the smallest non-overlined part, we can get
\begin{equation}\label{gen-G-N-O-new}
\sum_{\pi\in\mathcal{G}_O^{N}}q^{|\pi|}=\frac{1}{(q;q)_\infty}\sum_{n\geq 1}q^n(q^2;q^2)_{n-1}=(-q;q)_\infty\sum_{n\geq0} \frac{q^{2n+1}}{1-q^{2n+1}}\frac{1}{(q^2;q^2)_n},
\end{equation}
and
\begin{equation}\label{gen-E-N-O-new}
\sum_{\pi\in\mathcal{E}_O^{N}}q^{|\pi|}=(-q;q)_\infty\sum_{n\geq0} \frac{q^{2n+1}}{1-q^{2n+1}}\frac{1}{(q^2;q^2)_n}+\frac{1}{(q;q)_\infty}-(-q;q)_\infty.
\end{equation}
As a corollary of \eqref{gen-G-N-O-new}, we can get
\begin{core}
For $n\geq 1$, the number of overpartition of $n$ in  $\mathcal{G}_O^{N}$ equals the number of overpartitions of $n$ such that the largest non-overlined part is odd and the non-overlined parts less than   the largest non-overlined part are even.
\end{core}

In \cite{Cohen-1988}, Cohen observed the following identity:
\begin{equation*}
\sum_{n\geq 1}q^n(q^2;q^2)_{n-1}=\sum_{n\geq 1}\frac{(-1)^{n-1}q^{n^2}}{(q;q^2)_n}.
\end{equation*}
Combining with \eqref{gen-G-N-O-new}, we can get
\[\sum_{\pi\in\mathcal{G}_O^{N}}q^{|\pi|}=\frac{1}{(q;q)_\infty}\sum_{n\geq 1}\frac{(-1)^{n-1}q^{n^2}}{(q;q^2)_n}.\]

Inspired by the minimal excludant of a partition, Chern \cite{Chern-2021} investigated the maximal excludant of a partition. The maximal excludant of a partition $\pi$ is the largest nonnegative integer smaller than the largest part of $\pi$ that does not appear in $\pi$, denoted $maex(\pi)$. For $n\geq 1$, Chern defined
\[\sigma maex(n)=\sum_{\pi\in\mathcal{P}(n)}maex(\pi).\]
Let $\sigma L(n)$ be the sum of largest parts over all partitions of $n$, that is,
\[\sigma L(n)=\sum_{\pi\in\mathcal{P}(n)}LN(\pi).\]
Chern \cite{Chern-2021} obtained the following identity:
\[\sum_{n\geq 1}\left(\sigma L(n)-\sigma maex(n)\right)q^n=\frac{1}{(q;q)_\infty}\sum_{n\geq 1}q^n(q^2;q^2)_{n-1}.\]
Combining with \eqref{gen-G-N-O-new}, we have
\begin{core}\label{equiv-g-sigma}
For $n\geq 1$, the number of overpartitions of $n$ in $\mathcal{G}_O^{N}$ is  $\sigma L(n)-\sigma maex(n)$.
\end{core}

For $n\geq 1$, let ${E}_{NO}(n)$ denote the number of overpartitions of $n$ such that the size of the smallest non-overlined part equals the size of the largest overlined part, that is, ${E}_{NO}(n)$ is the number of overpartitions of $n$ with $SN(\pi)=LO(\pi)$.
For example, there are seven overpartitions counted by
${E}_{NO}(6)$, and so we have ${E}_{NO}(6)=7$.
\[(4,\overline{1},1),(\overline{3},3),(3,\overline{1},1,1),(\overline{2},2,2),(2,2,\overline{1},1),(2,\overline{1},1,1,1),(\overline{1},1,1,1,1,1).\]

It follows from \eqref{gen-G-N-O-new} and \eqref{gen-E-N-O-new} that
\begin{equation}\label{gen-E-N=O}
\sum_{n\geq 1}{E}_{NO}(n)q^{n}=\frac{1}{(q;q)_\infty}-(-q;q)_\infty.
\end{equation}
 Recall that the right hand side of \eqref{gen-E-N=O} is the generating function of $R(n)$, so we have
\begin{core}\label{core-NO=R}
For $n\geq1$,
\[{E}_{NO}(n)=R(n).\]
\end{core}

Let $D(n)$ denote the number of distinct partitions of $n$.
\begin{thm}\label{N-O-THM-1}
For $n\geq 1$,
\[{A}_{O\geq N}(n)-{B}_{O\geq N}(n)=D(n).\]
\end{thm}
It implies that
\begin{core}
For $n\geq 1$,
\[{A}_{O\geq N}(n)-{B}_{O\geq N}(n)>0.\]
\end{core}

Let ${H}'_{ON}(n)$ be the number of overpartitions $\pi$ of $n$ in $\mathcal{E}_O^{N}$ with $LN(\pi)=SN(\pi)$.

\begin{thm}\label{thm-e-o-gen-1}
For $n\geq 1$,
\[{A}_{O>N}(n)-{B}_{O>N}(n)={H}'_{ON}(n).\]
\end{thm}

It yields that
\begin{core}For $n\geq 1$,
\[{A}_{O>N}(n)-{B}_{O>N}(n)>0.\]
\end{core}

In Section 5, we will give proofs of \eqref{gen-G-N-O}, \eqref{gen-E-N-O}, \eqref{gen-G-N-O-new},  \eqref{gen-E-N-O-new}, Theorems \ref{N-O-THM-1} and \ref{thm-e-o-gen-1}, and we will present an analytic proof of \eqref{gen-E-N=O} and bijective proofs of Corollary \ref{equiv-g-sigma} and \ref{core-NO=R}.

\section{Preliminaries}
In this section, we collect some identities needed in this article. We also build an involution on $\overline{\mathcal{P}}(n)$, which plays a crucial in the proofs of Theorems \ref{EG-OG-THM}, \ref{EE-OE-THM}, \ref{N-O-THM-1} and \ref{thm-e-o-gen-1}.

%{\noindent \bf The $q$-binomial theorem \cite[Theorem 2.1]{Andrews-1976}:}
%\begin{equation}\label{chang-1}
%\sum_{n\geq0} \frac{(a;q)_{n}}{(q;q)_{n}}t^{n}=\frac{{(at;q)}_{\infty}}{(t;q)_{\infty}}.
%\end{equation}

{\noindent \bf Three identities due to Euler \cite{Euler-1748} (see also \cite[(1.2.5) and Corollary 2.2]{Andrews-1976}):}
\begin{equation}\label{Euler-1}
\sum_{n\geq0}\frac{t^n}{(q;q)_n}=\frac{1}{(t;q)_\infty},
\end{equation}
\begin{equation}\label{Euler-2}
\sum_{n\geq0}\frac{t^nq^{{n}\choose 2}}{(q;q)_n}=(-t;q)_\infty,
\end{equation}
and
\begin{equation}\label{Euler-new}
\frac{1}{(q;q^2)_\infty}=(-q;q)_\infty.
\end{equation}

Letting $q\rightarrow q^2$ and $t=q$ in \eqref{Euler-1}, we have
\begin{equation}\label{chang-1-1}
\sum_{n\geq0}\frac{q^n}{(q^2;q^2)_n}=\frac{1}{(q;q^2)_\infty}.
\end{equation}

Letting $q\rightarrow q^2$ and $t=q^2$ in \eqref{Euler-1}, we have
\begin{equation}\label{chang-1-2}
\sum_{n\geq0}\frac{q^{2n}}{(q^2;q^2)_n}=\frac{1}{(q^2;q^2)_\infty}.
\end{equation}

{\noindent \bf A formula due to Gauss \cite{Gauss-1866} (see also \cite[(2.2.13)]{Andrews-1976}):}
\begin{equation}\label{Gauss}
\sum_{n\geq 0}q^{\binom{n+1}{2}}=\frac{(q^2;q^2)_\infty}{(q;q^2)_\infty}.
\end{equation}

{\noindent \bf Heine's transformation formula \cite{Heine-1847} (see also \cite[Corollary 2.3]{Andrews-1976}):}
\begin{equation}\label{Heine}
\sum_{n\geq 0}\frac{(a;q)_n(b;q)_n}{(q;q)_n(c;q)_n}t^n=\frac{(b;q)_\infty(at;q)_\infty}{(c;q)_\infty(t;q)_\infty}\sum_{n\geq 0}\frac{(c/b;q)_n(t;q)_n}{(q;q)_n(at;q)_n}b^n.
\end{equation}

Let $\mathcal{D}$ denote the set of all distinct partitions. Then, we have
\begin{equation}\label{gen-distinct}
\sum_{\pi\in\mathcal{D}}t^{\ell(\pi)}q^{|\pi|}=(-tq;q)_\infty=1+\sum_{n\geq 1}tq^n(-tq;q)_{n-1}.
\end{equation}

We also need the following lemma.
\begin{lem}\label{pri-lem-1}
For $k\geq 0$,
\begin{equation}\label{pri-lem-eqn-1}
\sum_{m\geq k+1}\frac{q^m}{(q;q)_m}=\frac{1}{(q;q)_\infty}-\frac{1}{(q;q)_{k}}.
\end{equation}
\end{lem}
\pf It is clear that both sides of \eqref{pri-lem-eqn-1} are the generating function for the partitions with at least $k+1$ parts. This completes the proof.  \qed

We conclude this section with an involution $\mathcal{I}$ on $\overline{\mathcal{P}}(n)$.
\begin{defi}\label{defi-involution}
For $n\geq 1$, let $\pi$ be an overpartition  in $\overline{\mathcal{P}}(n)$. The map $\mathcal{I}$
is defined as follows:
\begin{itemize}
\item[{\rm(1)}] if $LN(\pi)>LO(\pi)$, then change the largest non-overlined part of $\pi$ to an overlined part{\rm;}

    \item[{\rm(2)}] if $LO(\pi)\geq LN(\pi)$, then change the largest overlined part of $\pi$ to a non-overlined part{\rm.}
\end{itemize}
\end{defi}

\section{Proofs of the results in Section 2.1}

In this section, we aim to show the results in Section 2.1.
We will prove \eqref{gen-G-O-N} in Section 4.1,  Corollary \ref{core-g-d2} and \ref{core-g-sigma} in Section 4.2, \eqref{gen-E-O-N} in Section 4.3, \eqref{gen-E-O=N} and Corollary \ref{core-ON=R} in Section 4.4, Theorem \ref{EG-OG-THM} in Section 4.5, and Theorem \ref{EE-OE-THM} in Section 4.6.

\subsection{Proofs of \eqref{gen-G-O-N}}

In this subsection, we give three proofs of \eqref{gen-G-O-N} in terms of separable overpartition classes, the largest non-overlined part and the smallest overlined part.

{\noindent \bf The first proof of \eqref{gen-G-O-N}.} For $m\geq 1$, define
\[\mathcal{BG}^O_N(m)=\left\{(1^m),(1^{m-1},\overline 2),\ldots,(1,\overline 2,\overline 3,\ldots,\overline{m-1},\overline m),(\overline 1,\overline 2,\overline 3,\ldots,\overline{m-1},\overline m)\right\}.\]
Set
  \[\mathcal{BG}^O_N=\bigcup_{m\geq 1}\mathcal{BG}^O_N(m).\]
Then, it is clear that $\mathcal{G}^O_N$ is a separable overpartition class and $\mathcal{BG}^O_N$ is the basis of $\mathcal{G}^O_N$. So, we get
\begin{align*}
\sum_{\pi\in\mathcal{G}_N^{O}}q^{|\pi|}&=1+\sum_{m\geq1}\frac{1}{(q;q)_m}\left[q^m+\sum_{k=1}^{m-1}q^{(m-k)+2+3+\cdots+(k+1)}+q^{\binom{m+1}{2}}\right]\\
&=1+\sum_{m\geq1}\frac{1}{(q;q)_m}\left[\sum_{k=0}^{m-1}q^{m+\binom{k+1}{2}}+q^{\binom{m+1}{2}}\right]\\
&=1+\sum_{m\geq1}\frac{q^m}{(q;q)_m}\sum_{k=0}^{m-1}q^{\binom{k+1}{2}}+\sum_{m\geq1}\frac{q^{\binom{m+1}{2}}}{(q;q)_m}\\
&=\sum_{k\geq0}q^{\binom{k+1}{2}}\sum_{m\geq k+1}\frac{q^m}{(q;q)_m}+\sum_{m\geq0}\frac{q^{\binom{m+1}{2}}}{(q;q)_m}.
\end{align*}
Combining with \eqref{Euler-new}, \eqref{Gauss} and Lemma \ref{pri-lem-1}, we have
\begin{align*}
\sum_{\pi\in\mathcal{G}_N^{O}}q^{|\pi|}&=\sum_{k\geq0}q^{\binom{k+1}{2}}\left[\frac{1}{(q;q)_\infty}-\frac{1}{(q;q)_{k}}\right]+\sum_{m\geq0}\frac{q^{\binom{m+1}{2}}}{(q;q)_m}\\
&=\frac{1}{(q;q)_\infty}\sum_{k\geq0}q^{\binom{k+1}{2}}-\sum_{k\geq0}\frac{q^{\binom{k+1}{2}}}{(q;q)_k}+\sum_{m\geq0}\frac{q^{\binom{m+1}{2}}}{(q;q)_m}\\
&=\frac{1}{(q;q)_\infty}\sum_{k\geq0}q^{\binom{k+1}{2}}\\
&=\frac{1}{(q;q)_\infty}\frac{(q^2;q^2)_\infty}{(q;q^2)_\infty}\\
&=(-q;q)_\infty^2.
\end{align*}
The proof is complete.  \qed

{\noindent \bf The second proof of \eqref{gen-G-O-N}.}
By considering the largest non-overlined part, we have
\[\sum_{\pi\in\mathcal{G}_N^{O}}q^{|\pi|}=\sum_{n\geq0}\frac{q^n}{(q;q)_n}(-q^{n+1};q)_\infty=(-q;q)_\infty\sum_{n\geq0}\frac{q^n}{(q;q)_n(-q;q)_n}=(-q;q)_\infty\sum_{n\geq0}\frac{q^n}{(q^2;q^2)_n}.\]
Combining with \eqref{chang-1-1} and \eqref{Euler-new}, we get
\[\sum_{\pi\in\mathcal{G}_N^{O}}q^{|\pi|}=(-q;q)_\infty\frac{1}{(q;q^2)_\infty}=(-q;q)_\infty^2.\]
The proof is complete. \qed

{\noindent \bf The third proof of \eqref{gen-G-O-N}.}
By virtue of the smallest overlined part, we get
\begin{align*}
\sum_{\pi\in\mathcal{G}_N^{O}}q^{|\pi|}&=\frac{1}{(q;q)_\infty}+\sum_{n\geq1}\frac{q^n(-q^{n+1};q)_\infty}{(q;q)_{n-1}}\\
&=\frac{1}{(q;q)_\infty}+(-q;q)_\infty\sum_{n\geq1}\frac{q^n(1-q^n)}{(q;q)_n(-q;q)_n}\\
&=\frac{1}{(q;q)_\infty}+(-q;q)_\infty\left[\sum_{n\geq0}\frac{q^n}{(q^2;q^2)_n}-\sum_{n\geq0}\frac{q^{2n}}{(q^2;q^2)_n}\right].
\end{align*}
Using \eqref{chang-1-1}, \eqref{chang-1-2} and \eqref{Euler-new}, we have
\[\sum_{\pi\in\mathcal{G}_N^{O}}q^{|\pi|}=\frac{1}{(q;q)_\infty}+(-q;q)_\infty\left[\frac{1}{(q;q^2)_\infty}-\frac{1}{(q^2;q^2)_\infty}\right]=(-q;q)_\infty^2.\]
This completes the proof.  \qed

\subsection{Combinatorial proofs of Corollary \ref{core-g-d2} and \ref{core-g-sigma}}

This subsection is devoted to giving combinatorial proofs of Corollary \ref{core-g-d2} and \ref{core-g-sigma}. We first prove Corollary \ref{core-g-d2} by a bijective method, which can be seen as a combinatorial proof of \eqref{gen-G-O-N}. For $n,k\geq 0$, let $\mathcal{G}_N^{O}(n,k)$ be the set of overpartitions of $n$ in $\mathcal{G}_N^{O}$ with exactly $k$ overlined parts and let $\mathcal{D}_{2}(n,k)$ be the set of pairs $(\alpha,\beta)$ counted by $D_2(n)$ with $\ell(\alpha)-\ell(\beta)=k$ or $-k-1$. To prove Corollary \ref{core-g-d2}, it suffices to show the following theorem.
\begin{thm}\label{bijective-ab}
For $n,k\geq 0$, there is a bijection between $\mathcal{G}_N^{O}(n,k)$ and $\mathcal{D}_{2}(n,k)$.
\end{thm}
We will give a combinatorial proof of Theorem \ref{bijective-ab} by a modification of Sylvester's bijective proof of Jacobi's product identity \cite{Franklin-Sylvester-1882}, which was later rediscovered by Wright \cite{Wright-1965}.

{\noindent \bf Proof of Theorem \ref{bijective-ab}.} Let $\pi=(\pi_1,\pi_2,\ldots,\pi_m)$ be an overpartition in $\mathcal{G}_N^{O}(n,k)$. Under the definition of $\mathcal{G}_N^{O}(n,k)$, we see that $\pi_1,\ldots,\pi_k$ are overlined parts, $\pi_{k+1},\ldots,\pi_m$ are non-overlined parts and
\[|\pi_1|>\cdots>|\pi_k|>|\pi_{k+1}|\geq\cdots\geq|\pi_m|.\]
Set $\lambda_i=|\pi_i|-(k-i+1)$ for $1\leq i\leq k$, and  $\lambda_i=\pi_i$ for $k+1\leq i\leq m$. It is clear that $\lambda=(\lambda_1,\lambda_2,\ldots,\lambda_m)$ is a partition of $n-{{k+1}\choose 2}$. Let $\lambda'=(\lambda'_1,\lambda'_2,\ldots)$ be the conjugate of $\lambda$. Set $s$ be the largest integer such that $\lambda_{s}\geq k+s$ with the convention that $\lambda_0=+\infty$. By the choice of $s$, we have
\begin{equation}\label{conjugate}
\lambda_{s+1}\leq k+s,\ \lambda'_{k+s}\geq s\text{ and }\lambda'_{k+s+1}\leq s.
\end{equation}
For $i\geq 1$, we set $\mu_i=\lambda'_i+k-i+1$ and $\nu_i=\lambda_i-k-i$. Clearly, we have $\mu_1>\mu_2>\cdots$ and  $\nu_1>\nu_2>\cdots$. Appealing to \eqref{conjugate}, we derive that
\[\mu_{k+s}=\lambda'_{k+s}+k-(k+s)+1\geq s+k-(k+s)+1=1,\]
\[\mu_{k+s+1}=\lambda'_{k+s+1}+k-(k+s+1)+1\leq s+k-(k+s+1)+1=0,\]
and
\[\nu_{s+1}=\lambda_{s+1}-k-(s+1)\leq (k+s)-k-(s+1)=-1.\]
Setting $\mu=(\mu_1,\mu_2,\ldots,\mu_{k+s})$, we see that $\mu$ is a distinct partition. Note that $\lambda_{s}\geq k+s$, then we consider the following two cases.

Case 1: $\lambda_s\geq k+s+1$. In such case, we have $\nu_{s}=\lambda_s-k-s\geq (k+s+1)-k-s=1$. Then, we set $\beta=(\nu_1,\nu_2,\ldots,\nu_s)$, which is a distinct partition. Setting $\alpha=\mu$, then $(\alpha,\beta)$ is a pair of distinct partitions with $\ell(\alpha)-\ell(\beta)=k$.

Case 2: $\lambda_s=k+s$. In such case, we have $\nu_{s}=\lambda_s-k-s=(k+s)-k-s=0$. Then, we set $\alpha=(\nu_1,\nu_2,\ldots,\nu_{s-1})$, which a distinct partition. Setting $\beta=\mu$, then $(\alpha,\beta)$ is a pair of distinct partitions with $\ell(\alpha)-\ell(\beta)=-k-1$.

 In either case, we get a pair of distinct partitions $(\alpha,\beta)$ with $\ell(\alpha)-\ell(\beta)=k$ or $-k-1$. It is clear that $|\alpha|+|\beta|=|\lambda|+{{k+1}\choose 2}=n$. So, we have $(\alpha,\beta)\in\mathcal{D}_{2}(n,k)$. Obviously, the process above is reversible. The proof is complete. \qed

For example, let $\pi=(\overline{10},\overline{8},\overline{7},6,4,4,2,1)$ be an overpartition in $\mathcal{G}_N^{O}(42,3)$. Set $\lambda_1=|\overline{10}|-3=7$, $\lambda_2=|\overline{8}|-2=6$, $\lambda_3=|\overline{7}|-1=6$, $\lambda_4=6$, $\lambda_5=4$, $\lambda_6=4$, $\lambda_7=2$ and $\lambda_8=1$. Then, we get $\lambda=(7,6,6,6,4,4,2,1)$, which is a partition of $36$. The conjugate of $\lambda$ is $\lambda'=(8,7,6,6,4,4,1)$. It can be check that $s=3$ is the largest integer such that $\lambda_{s}=\lambda_3=6\geq 3+s=6$. Moreover, we have $\lambda_{s}=\lambda_3=6=3+s$. Then, we can get
\[\alpha=(3,1)\text{ and }\beta=(11,9,7,6,3,2),\]
where $\alpha_i=\lambda_i-3-i$ for $1\leq i\leq 2$, and $\beta_i=\lambda'_i+3-i+1$ for $1\leq i\leq 6$.
Clearly, $(\alpha,\beta)$ is a pair of distinct partitions in $\mathcal{D}_{2}(42,3)$. The same process to get $(\alpha,\beta)$ could be run in reverse.

For another example, let $\pi=(7,6,6,6,4,4,2,1)$ be an overpartition in $\mathcal{G}_N^{O}(36,0)$. Note that $k=0$, then we have $\lambda=\pi=(7,6,6,6,4,4,2,1)$. The conjugate of $\lambda$ is $\lambda'=(8,7,6,6,4,4,1)$. It can be check that $s=4$ is the largest integer such that $\lambda_{s}=\lambda_4=6\geq 0+s=4$. Moreover, we have $\lambda_{s}=\lambda_4=6\geq 0+s+1=5$. Then, we can get
\[\alpha=(8,6,4,3)\text{ and }\beta=(6,4,3,2),\]
where $\alpha_i=\lambda'_i+0-i+1$ for $1\leq i\leq 4$, and $\beta_i=\lambda_i-0-i$ for $1\leq i\leq 4$. Clearly, $(\alpha,\beta)$ is a pair of distinct partitions in $\mathcal{D}_{2}(36,0)$. The same process to get $(\alpha,\beta)$ could be run in reverse.

Then, we proceed to give the combinatorial proof of Corollary \ref{core-g-sigma}.

{\noindent \bf Combinatorial proof of Corollary \ref{core-g-sigma}.} For $n,k\geq 0$, we define
\[\mathcal{M}(n,k)=\left\{\mu\mid|\mu|=n,0\leq k\leq mex(\mu)-1\right\}.\]
  We just need to build a bijection between $\mathcal{G}_N^{O}(n,k)$ and $\mathcal{M}(n,k)$.

Let $\pi=(\pi_1,\pi_2,\ldots,\pi_m)$ be an overpartition in $\mathcal{G}_N^{O}(n,k)$. With the same argument in the proof of Theorem \ref{bijective-ab}, we can get a partition $\lambda=(\lambda_1,\lambda_2,\ldots,\lambda_m)$ of $n-{{k+1}\choose 2}$, where \[\lambda_i=|\pi_i|-(k-i+1)\text{ for }1\leq i\leq k,\text{ and  }\lambda_i=\pi_i\text{ for }k+1\leq i\leq m.\]
 Then, we add $1,2,\ldots,k$ as parts into $\lambda$ and denote the resulting partition by  $\mu$.
 Namely,
 \[f_t(\mu)=f_t(\lambda)+1\text{ for }1\leq t\leq k,\text{ and }f_t(\mu)=f_t(\lambda)\text{ otherwise.}\]
 Clearly, $\mu$ is a partition of $n$ with $mex(\mu)\geq k+1$, and so $\mu$ is a partition  in $\mathcal{M}(n,k)$. Obviously, the process above is reversible. The proof is complete. \qed

\subsection{Proofs of \eqref{gen-E-O-N}}

With the similar arguments in the proofs of \eqref{gen-G-O-N} given in Section 4.1, we present three proofs of \eqref{gen-E-O-N} in this subsection.

{\noindent \bf The first proof of \eqref{gen-E-O-N}.} For $m\geq 1$, define
\[\mathcal{BE}^O_N(m)=\left\{(1^m), (1^{m-1},\overline{1}),\ldots,(1,\overline{1},\overline{2},\ldots,\overline{m-1},\overline{m}),(\overline 1,\overline 2,\overline 3,\ldots,\overline{m-1},\overline m)\right\}.\]
Set
  \[\mathcal{BE}^O_N=\bigcup_{m\geq 1}\mathcal{BE}^O_N(m).\]
Then, it is clear that $\mathcal{E}^O_N$ is a separable overpartition class and $\mathcal{BE}^O_N$ is the basis of $\mathcal{E}^O_N$. So, we get
\begin{align*}
\sum_{\pi\in\mathcal{E}_N^{O}}q^{|\pi|}&=1+\sum_{m\geq1}\frac{1}{(q;q)_m}\left[q^m+\sum_{k=1}^{m-1}q^{(m-k)+1+2+\cdots+k}+q^{\binom{m+1}{2}}\right]\\
&=1+\sum_{m\geq1}\frac{1}{(q;q)_m}\sum_{k=0}^{m}q^{m+\binom{k}{2}}\\
&=\sum_{k\geq0}q^{\binom{k}{2}}\sum_{m\geq k}\frac{q^m}{(q;q)_m}\\
&=\sum_{m\geq 0}\frac{q^m}{(q;q)_m}+\sum_{k\geq1}q^{\binom{k}{2}}\sum_{m\geq k}\frac{q^m}{(q;q)_m}\\
&=\frac{1}{(q;q)_\infty}+\sum_{k\geq1}q^{\binom{k}{2}}\left[\frac{1}{(q;q)_\infty}-\frac{1}{(q;q)_{k-1}}\right],
\end{align*}
where the final equation follows from \eqref{Euler-1} with $t=q$ and Lemma \ref{pri-lem-1}. Appealing to \eqref{Euler-2} with $t=q$, \eqref{Euler-new} and \eqref{Gauss}, we have
\begin{align*}
\sum_{\pi\in\mathcal{E}_N^{O}}q^{|\pi|}&=\frac{1}{(q;q)_\infty}+\frac{1}{(q;q)_\infty}\sum_{k\geq0}q^{\binom{k+1}{2}}-\sum_{k\geq0}\frac{q^{\binom{k+1}{2}}}{(q;q)_{k}}\\
&=\frac{1}{(q;q)_\infty}+\frac{1}{(q;q)_\infty}\frac{(q^2;q^2)_\infty}{(q;q^2)_\infty}-(-q;q)_\infty\\
&=\frac{1}{(q;q)_\infty}+(-q;q)_\infty^2-(-q;q)_\infty.
\end{align*}
This completes the proof.  \qed

{\noindent \bf The second proof of \eqref{gen-E-O-N}.}
In view of the largest non-overlined part, we have
\begin{align*}
\sum_{\pi\in\mathcal{E}_N^{O}}q^{|\pi|}&=(-q;q)_\infty+\sum_{n\geq1}\frac{q^n}{(q;q)_n}(-q^{n};q)_\infty\\
&=(-q;q)_\infty+(-q;q)_\infty\sum_{n\geq1}\frac{q^n(1+q^n)}{(q;q)_n(-q;q)_n}\\
&=(-q;q)_\infty+(-q;q)_\infty\left[\sum_{n\geq0}\frac{q^n}{(q^2;q^2)_n}+\sum_{n\geq0}\frac{q^{2n}}{(q^2;q^2)_n}-2\right].
\end{align*}
It follows from \eqref{chang-1-1}, \eqref{chang-1-2} and \eqref{Euler-new} that
\begin{align*}
\sum_{\pi\in\mathcal{E}_N^{O}}q^{|\pi|}&=(-q;q)_\infty+(-q;q)_\infty\left[\frac{1}{(q;q^2)_\infty}+\frac{1}{(q^2;q^2)_\infty}-2\right]\\
&=(-q;q)_\infty^2+\frac{1}{(q;q)_\infty}-(-q;q)_\infty.
\end{align*}
The proof is complete.  \qed

{\noindent \bf The third proof of \eqref{gen-E-O-N}.}
In light of the smallest overlined part, we have
\begin{align*}
\sum_{\pi\in\mathcal{E}_N^{O}}q^{|\pi|}&=\frac{1}{(q;q)_\infty}+\sum_{n\geq1}q^n(-q^{n+1};q)_\infty\frac{1}{(q;q)_{n}}\\
&=\frac{1}{(q;q)_\infty}+(-q;q)_\infty\sum_{n\geq 1}\frac{q^n}{(-q;q)_{n}(q;q)_{n}}\\
&=\frac{1}{(q;q)_\infty}+(-q;q)_\infty\left[\sum_{n\geq0}\frac{q^n}{(q^2;q^2)_n}-1\right].
\end{align*}
Combining with \eqref{chang-1-1} and \eqref{Euler-new}, we arrive at
\begin{align*}
\sum_{\pi\in\mathcal{E}_N^{O}}q^{|\pi|}&=\frac{1}{(q;q)_\infty}+(-q;q)_\infty\left[\frac{1}{(q;q^2)_\infty}-1\right]\\
&=\frac{1}{(q;q)_\infty}+(-q;q)_\infty^2-(-q;q)_\infty.
\end{align*}
This completes the proof.  \qed

\subsection{Proofs of \eqref{gen-E-O=N} and Corollary \ref{core-ON=R}}

In this subsection, we first give an analytic proof of \eqref{gen-E-O=N} and then we give a combinatorial proof of Corollary \ref{core-ON=R}.

{\noindent \bf Analytic proof of \eqref{gen-E-O=N}.} By considering the size of the smallest overlined part and the largest non-overlined part, we can get
\begin{align*}
\sum_{n\geq 1}{E}_{ON}(n)q^{n}&=\sum_{n\geq1}{q^{n}}(-q^{n+1};q)_\infty\frac{q^{n}}{(q;q)_n}\\
&=(-q;q)_\infty\sum_{n\geq1}\frac{q^{2n}}{(-q;q)_{n}(q;q)_{n}}\\
&=(-q;q)_\infty\left[\sum_{n\geq 0}\frac{q^{2n}}{(q^2;q^2)_n}-1\right].
\end{align*}
Combining with \eqref{chang-1-2}, we have
\[
\sum_{n\geq 1}{E}_{ON}(n)q^{n}=(-q;q)_\infty\left[\frac{1}{(q^2;q^2)_\infty}-1\right]=\frac{1}{(q;q)_\infty}-(-q;q)_\infty.\]
The proof is complete.  \qed

Now, we proceed to give a combinatorial proof of Corollary \ref{core-ON=R}, which can be seen as a bijective proof of \eqref{gen-E-O=N}.

{\noindent \bf Combinatorial proof of Corollary \ref{core-ON=R}.} Let $\pi$ be an overpartition counted by ${E}_{ON}(n)$. We define $\lambda=\phi_{ON}(\pi)$ by changing the overlined parts in $\pi$ to non-overlined parts. Clearly, $\lambda$ is a partition enumerated by $R(n)$.

Conversely, let $\lambda$ be a partition enumerated by $R(n)$. Assume that $k$ is the largest integer such that $k$ appears at least twice in $\lambda$. Then, we define $\pi=\psi_{ON}(\lambda)$ by changing a $k$ in $\lambda$ to $\overline k$ and changing the parts greater than $k$ to overlined parts. It is obvious that $\pi$ is an overpartition counted by ${E}_{ON}(n)$. Clearly, $\psi_{ON}$ is the inverse map of $\phi_{ON}$. Thus, we complete the proof.  \qed

\subsection{Proof of Theorem \ref{EG-OG-THM}}

This subsection is devoted to showing Theorem \ref{EG-OG-THM}. Let $p^{e}_o(n)$ (resp. $p^{o}_e(n)$) be the number of partitions of $n$ such that the largest part is even (resp. odd) and the smallest part is odd (resp. even).
Then, it suffices to show that
\begin{equation}\label{lem-excess-1}
\sum_{n\geq 0}\left({A}_{N\geq O}(n)-{B}_{N\geq O}(n)\right)q^n=1+2\sum_{n\geq 1}\left(p^{e}_o(n)-p^{o}_e(n)\right)q^n,
\end{equation}
and for $n\geq 1$,
\begin{equation}\label{lem-excess-2}
p^{e}_o(n)-p^{o}_e(n)\geq 0\text{ with strict inequality if }n\geq 3.
\end{equation}

We first give an analytic proof and a combinatorial proof of \eqref{lem-excess-1}.

{\noindent \bf Analytic proof of \eqref{lem-excess-1}.} In virtue of the smallest overlined part, we can get
\begin{align}
\sum_{n\geq 0}\left({A}_{N\geq O}(n)-{B}_{N\geq O}(n)\right)q^n&=\frac{1}{(-q;q)_\infty}+\sum_{n\geq1} q^n(-q^{n+1};q)_\infty\frac{1}{(q;q)_{n-1}}\frac{1}{(-q^n;q)_\infty}\nonumber\\
&=\frac{1}{(-q;q)_\infty}+\sum_{n\geq1}\frac{1}{(q;q)_{n-1}}\frac{q^n}{1+q^n}.\label{proof-lem-excess-1}
\end{align}

Utilizing \eqref{Euler-1} with $t=-q$, we have
\begin{equation}\label{proof-lem-excess-2}
\frac{1}{(-q;q)_\infty}=\sum_{n\geq0}\frac{(-q)^n}{(q;q)_n}=1+\sum_{n\geq1}\frac{q^{2n}}{(q;q)_{2n}}-\sum_{n\geq0}\frac{q^{2n+1}}{(q;q)_{2n+1}}.
\end{equation}
Clearly,
\[\sum_{n\geq1}\frac{q^{2n}}{(q;q)_{2n}}\ \left(\text{resp. }\sum_{n\geq0}\frac{q^{2n+1}}{(q;q)_{2n+1}}\right)\]
is the generating function for the nonempty partitions such that  the largest part is even (resp. odd).

Setting $t=q^{m+1}$ in \eqref{Euler-1}, we have
\[\sum_{n\geq0}\frac{q^{(m+1)n}}{(q;q)_n}=\frac{1}{(q^{m+1};q)_\infty}.\]
So, we get
\begin{align}
\sum_{n\geq1}\frac{1}{(q;q)_{n-1}}\frac{q^n}{1+q^n}&=\sum_{n\geq0}\frac{1}{(q;q)_n}\frac{q^{n+1}}{1+q^{n+1}}\nonumber\\
&=\sum_{n\geq 0}\frac{1}{(q;q)_n}\sum_{m\geq 0}(-1)^mq^{(n+1)(m+1)}\nonumber\\
&=\sum_{m\geq0}(-1)^mq^{m+1}\sum_{n\geq 0}\frac{q^{(m+1)n}}{(q;q)_n}\nonumber\\
&=\sum_{m\geq 0}(-1)^mq^{m+1}\frac{1}{(q^{m+1};q)_\infty}\nonumber\\
&=\sum_{m\geq 0}\frac{q^{2m+1}}{(q^{2m+1};q)_\infty}-\sum_{m\geq0}\frac{q^{2m+2}}{(q^{2m+2};q)_\infty}.\label{proof-lem-excess-3}
\end{align}
Clearly,
\[\sum_{m\geq 0}\frac{q^{2m+1}}{(q^{2m+1};q)_\infty}\ \left(\text{resp. }\sum_{m\geq0}\frac{q^{2m+2}}{(q^{2m+2};q)_\infty}\right)\]
is the generating function for the nonempty partitions such that  the smallest part is odd (resp. even).

Then, we obtain that
\[\sum_{n\geq1}\frac{q^{2n}}{(q;q)_{2n}}-\sum_{m\geq0}\frac{q^{2m+2}}{(q^{2m+2};q)_\infty}\]
and
\[\sum_{m\geq 0}\frac{q^{2m+1}}{(q^{2m+1};q)_\infty}-\sum_{n\geq0}\frac{q^{2n+1}}{(q;q)_{2n+1}}\]
are both equal to the generating function for the nonempty partitions such that the largest part is even  and the smallest part is odd  subtract the generating function for the non-empty partitions such that the largest part is odd and the smallest part is even, that is,
\begin{align}
\sum_{n\geq 1}\left(p^{e}_o(n)-p^{o}_e(n)\right)q^n&=\sum_{n\geq1}\frac{q^{2n}}{(q;q)_{2n}}-\sum_{m\geq0}\frac{q^{2m+2}}{(q^{2m+2};q)_\infty}\label{proof-lem-excess-4}\\
&=\sum_{m\geq 0}\frac{q^{2m+1}}{(q^{2m+1};q)_\infty}-\sum_{n\geq0}\frac{q^{2n+1}}{(q;q)_{2n+1}}.\label{proof-lem-excess-5}
\end{align}

Combining \eqref{proof-lem-excess-1}, \eqref{proof-lem-excess-2}, \eqref{proof-lem-excess-3}, \eqref{proof-lem-excess-4} and \eqref{proof-lem-excess-5}, we arrive at  \eqref{lem-excess-1}. This completes the proof.  \qed

Before giving the combinatorial proof of \eqref{lem-excess-1}, we introduce the following notations, which will be also used in the  combinatorial proof of Theorem  \ref{EE-OE-THM}.

For $n\geq 1$, we introduce the following three sets.
\begin{itemize}
\item Let $\mathcal{C}_{NO}(n)$ be the set of overpartitions $\pi$ of $n$ such that $LN(\pi)>LO(\pi)\geq 1$.

\item Let $\mathcal{F}_{NO}(n)$ be the set of overpartitions $\pi$ of $n$ such that there are at least two overlined parts in $\pi$ and $LO(\pi)\geq LN(\pi)$.

\item Let $\mathcal{H}_{NO}(n)$ be the set of overpartitions $\pi$ of $n$ such that there is exactly one overlined part in $\pi$ and $LO(\pi)\geq LN(\pi)$.

\end{itemize}
Clearly, we have
\[\overline{\mathcal{P}}(n)=\mathcal{C}_{NO}(n)\bigcup\mathcal{F}_{NO}(n)\bigcup\mathcal{H}_{NO}(n)\bigcup\mathcal{P}(n).\]

By restricting the involution $\mathcal{I}$ defined in Definition \ref{defi-involution} on $\mathcal{C}_{NO}(n)\bigcup\mathcal{F}_{NO}(n)$, it is easy to get the following lemma.
\begin{lem}\label{lem-involution-CF}
For $n\geq 1$, the map $\mathcal{I}$ is a bijection between $\mathcal{C}_{NO}(n)$ and $\mathcal{F}_{NO}(n)$. Moreover,  for an overpartition $\pi \in \mathcal{C}_{NO}(n)$, let $\lambda=\mathcal{I}(\pi)$. We have \[\ell_{N\geq O}(\lambda)=\ell_{N\geq O}(\pi)-1\text{ and }\ell_{N>O}(\lambda)=\ell_{N>O}(\pi)-1.\]
\end{lem}

Then, we proceed to give the combinatorial proof of \eqref{lem-excess-1}.

{\noindent \bf Combinatorial proof of \eqref{lem-excess-1}.}  Note that
\begin{equation*}
\sum_{n\geq 0}\left({A}_{N\geq O}(n)-{B}_{N\geq O}(n)\right)q^n=1+\sum_{n\geq 1}\sum_{\pi\in\overline{\mathcal{P}}(n)}(-1)^{\ell_{N\geq O}(\pi)}q^{n},
\end{equation*}
it is equivalent to showing that for $n\geq 1$,
 \begin{equation}\label{lem-excess-1-combinatorial-1}
 \sum_{\pi\in\overline{\mathcal{P}}(n)}(-1)^{\ell_{N\geq O}(\pi)}=2\left(p^{e}_o(n)-p^{o}_e(n)\right).
\end{equation}
To do this, we are required to compute  \[\sum_{\pi\in\mathcal{C}_{NO}(n)\bigcup\mathcal{F}_{NO}(n)}(-1)^{\ell_{N\geq O}(\pi)},\ \sum_{\pi\in\mathcal{H}_{NO}(n)}(-1)^{\ell_{N\geq O}(\pi)}\text{ and }\sum_{\pi\in\mathcal{P}(n)}(-1)^{\ell_{N\geq O}(\pi)}.\]

By Lemma \ref{lem-involution-CF}, we have
\begin{equation}\label{lem-excess-1-combinatorial-2}
 \sum_{\pi\in\mathcal{C}_{NO}(n)\bigcup\mathcal{F}_{NO}(n)}(-1)^{\ell_{N\geq O}(\pi)}=0.
\end{equation}

Then, we divide $\mathcal{H}_{NO}(n)$ into four disjoint subsets.
\begin{itemize}
\item Let $\mathcal{H}_{NO}^{(1)}(n)$ (resp. $\mathcal{H}^{(2)}_{NO}(n)$) be the set of overpartitions $\pi$ in $\mathcal{H}_{NO}(n)$ such that $\ell(\pi)$ is even and there is an even (resp. odd) number of non-overlined parts of size $LO(\pi)$ in $\pi$.

\item Let $\mathcal{H}_{NO}^{(3)}(n)$ (resp. $\mathcal{H}^{(4)}_{NO}(n)$) be the set of overpartitions $\pi$ in $\mathcal{H}_{NO}(n)$ such that $\ell(\pi)$ is odd and there is an even (resp. odd) number of non-overlined parts of size $LO(\pi)$ in $\pi$.
\end{itemize}
Obviously, we have
\[\mathcal{H}_{NO}(n)=\mathcal{H}_{NO}^{(1)}(n)\bigcup\mathcal{H}_{NO}^{(2)}(n)\bigcup\mathcal{H}_{NO}^{(3)}(n)\bigcup\mathcal{H}_{NO}^{(4)}(n).\]

For an overpartition $\pi\in \mathcal{H}_{NO}(n)$, it is clear that $\ell_{N\geq O}(\pi)$ is the number of  non-overlined parts of size $LO(\pi)$ in $\pi$. So, we have $(-1)^{\ell_{N\geq O}(\pi)}=1$ if $\pi\in\mathcal{H}_{NO}^{(1)}(n)\bigcup\mathcal{H}_{NO}^{(3)}(n)$, and $(-1)^{\ell_{N\geq O}(\pi)}=-1$ if $\pi\in\mathcal{H}_{NO}^{(2)}(n)\bigcup\mathcal{H}_{NO}^{(4)}(n)$.

We also divide $\mathcal{P}(n)$ into four disjoint subsets.
\begin{itemize}
\item Let $\mathcal{P}^{(1)}(n)$ (resp. $\mathcal{P}^{(2)}(n)$) be the set of partitions $\pi$ of $n$ such that $\ell(\pi)$ is even and there is an odd (resp. even) number of parts of size $LN(\pi)$ in $\pi$.

\item Let $\mathcal{P}^{(3)}(n)$ (resp. $\mathcal{P}^{(4)}(n)$) be the set of partitions $\pi$ of $n$ such that $\ell(\pi)$ is odd and there is an odd (resp. even) number of parts of size $LN(\pi)$ in $\pi$.
\end{itemize}
Obviously, we have
\[\mathcal{P}(n)=\mathcal{P}^{(1)}(n)\bigcup\mathcal{P}^{(2)}(n)\bigcup\mathcal{P}^{(3)}(n)\bigcup\mathcal{P}^{(4)}(n).\]

For a partition $\pi\in \mathcal{P}(n)$, it is clear that $\ell_{N\geq O}(\pi)=\ell(\pi)$. So, we have $(-1)^{\ell_{N\geq O}(\pi)}=1$ if $\pi\in\mathcal{P}^{(1)}(n)\bigcup\mathcal{P}^{(2)}(n)$, and $(-1)^{\ell_{N\geq O}(\pi)}=-1$ if $\pi\in\mathcal{P}^{(3)}(n)\bigcup\mathcal{P}^{(4)}(n)$.

By restricting the involution $\mathcal{I}$ defined in Definition \ref{defi-involution} on $\mathcal{H}_{NO}(n)\bigcup\mathcal{P}(n)$, we find that the map $\mathcal{I}$ is a bijection between $\mathcal{H}_{NO}^{(i)}(n)$ and $\mathcal{P}^{(i)}(n)$ for $1\leq i\leq 4$. So, we get
\begin{equation}\label{lem-excess-1-combinatorial-3}
 \sum_{\pi\in\mathcal{H}_{NO}(n)\bigcup\mathcal{P}(n)}(-1)^{\ell_{N\geq O}(\pi)}=2\left(\sum_{\pi\in\mathcal{P}^{(1)}(n)}1-\sum_{\pi\in\mathcal{P}^{(4)}(n)}1\right).
\end{equation}

By considering the conjugate, we can get
\begin{equation*}\label{lem-excess-1-combinatorial-5}
\sum_{\pi\in\mathcal{P}^{(1)}(n)}1=p^{e}_o(n)\text{ and }\sum_{\pi\in\mathcal{P}^{(4)}(n)}1=p^{o}_e(n).
\end{equation*}
Combining with  \eqref{lem-excess-1-combinatorial-2} and \eqref{lem-excess-1-combinatorial-3}, we arrive at \eqref{lem-excess-1-combinatorial-1}, and thus the proof is complete.  \qed

We conclude this section with the proof of \eqref{lem-excess-2}.

{\noindent \bf Proof of \eqref{lem-excess-2}.} It can be checked that
\[p^{e}_o(1)=p^{o}_e(1)=p^{e}_o(2)=p^{o}_e(2)=0,\]
which implies that \eqref{lem-excess-2} holds for $n=1,2$.

For $n\geq 3$, let $\mathcal{P}^{o}_e(n)$ be the set of partitions counted by $p^{o}_e(n)$ and let $\hat{\mathcal{P}}^{e}_o(n)$ be set of partitions $\pi$ enumerated by $p^{e}_o(n)$ such that $\ell(\pi)\geq 4$, $\ell(\pi)$ is even, $\pi_{\frac{\ell(\pi)}{2}}$ is odd and $f_{1}(\pi)\geq\frac{\ell(\pi)}{2}$. Then, we establish a bijection $\phi^o_{e}$ between $\mathcal{P}^{o}_e(n)$ and $\hat{\mathcal{P}}^{e}_o(n)$.

For a partition $\lambda\in\mathcal{P}^{o}_e(n)$, we set
\[\pi_i=\lambda_i-1\text{ for }1\leq i\leq \ell(\lambda),\text{ and }\pi_i=1\text{ for }\ell(\lambda)+1\leq i\leq 2\ell(\lambda).\]
Then, we define $\phi^o_{e}(\lambda)=(\pi_1,\pi_2,\ldots,\pi_{2\ell(\lambda)})$, which is a partition in $\hat{\mathcal{P}}^{e}_o(n)$.

Conversely, For a partition $\pi\in\hat{\mathcal{P}}^{e}_o(n)$, we set
\[\lambda_i=\pi_i+1\text{ for }1\leq i\leq \frac{\ell(\pi)}{2}.\]
Then, we define $\psi^o_{e}(\pi)=(\lambda_1,\lambda_2,\ldots,\lambda_{\frac{\ell(\pi)}{2}})$, which is a partition in $\mathcal{P}^{o}_e(n)$. Clearly, $\psi^o_{e}$ is the inverse map of $\phi^o_{e}$.

Now, we have built the bijection $\phi^o_{e}$ between $\mathcal{P}^{o}_e(n)$ and $\hat{\mathcal{P}}^{e}_o(n)$. Note that the partitions in $\hat{\mathcal{P}}^{e}_o(n)$ are counted by  $p^{e}_o(n)$. In order to show that \eqref{lem-excess-2} holds for $n\geq 3$, it remains to find a partition which is counted by  $p^{e}_o(n)$ but is not in $\hat{\mathcal{P}}^{e}_o(n)$. For a partition $\pi\in\hat{\mathcal{P}}^{e}_o(n)$, by definition, we have $\ell(\pi)\geq 4$. So, we just need to find a partition counted by  $p^{e}_o(n)$ with less than four parts. There are two cases.

Case 1: $n$ is odd. In such case, $(n-1,1)$ is counted by  $p^{e}_o(n)$ but is not in $\hat{\mathcal{P}}^{e}_o(n)$.

Case 2: $n$ is even. In such case, $(n-2,1,1)$ is counted by  $p^{e}_o(n)$ but is not in $\hat{\mathcal{P}}^{e}_o(n)$.

Thus, the proof is complete.  \qed

\subsection{Proofs of Theorem  \ref{EE-OE-THM}}

In this subsection, we will give an analytic proof and a combinatorial proof of Theorem  \ref{EE-OE-THM}.
We first give an analytic proof of Theorem \ref{EE-OE-THM}.

{\noindent \bf Analytic proof of Theorem \ref{EE-OE-THM}.} By considering the smallest overlined part, we can get
\begin{align*}
\sum_{n\geq 0}\left({A}_{N>O}(n)-{B}_{N>O}(n)\right)q^n&=\frac{1}{(-q;q)_\infty}+\sum_{n\geq1} q^n(-q^{n+1};q)_\infty\frac{1}{(q;q)_{n}}\frac{1}{(-q^{n+1};q)_\infty}\\
&=\frac{1}{(-q;q)_\infty}+\sum_{n\geq 1}\frac{q^n}{(q;q)_{n}}\\
&=1+2\sum_{n\geq1}\frac{q^{2n}}{(q;q)_{2n}},
\end{align*}
where the final equation follows from \eqref{proof-lem-excess-2}. This completes the proof.  \qed

Then, we give a combinatorial proof of Theorem \ref{EE-OE-THM} with a similar argument in the combinatorial of \eqref{lem-excess-1}.

{\noindent \bf Combinatorial proof of Theorem \ref{EE-OE-THM}.}
We find that it is equivalent to showing that for $n\geq 1$,
 \begin{equation}\label{EE-OE-THM-combinatorial-0}
 \sum_{\pi\in\overline{\mathcal{P}}(n)}(-1)^{\ell_{N> O}(\pi)}=2p_e(n).
\end{equation}

For $n\geq1$, let  $\mathcal{P}_{e}(n)$ (resp. $\mathcal{P}_{o}(n)$) be the set of partitions of $n$ with an even (resp. odd) number of parts.
Obviously, we have
\[\mathcal{P}(n)=\mathcal{P}_{e}(n)\bigcup\mathcal{P}_{o}(n).\]

For a partition $\pi\in \mathcal{P}(n)$, it is clear that $\ell_{N> O}(\pi)=\ell(\pi)$. So, we have $(-1)^{\ell_{N> O}(\pi)}=1$ if $\pi\in\mathcal{P}_{e}(n)$, and $(-1)^{\ell_{N> O}(\pi)}=-1$ if $\pi\in\mathcal{P}_{o}(n)$.  So we get
\begin{equation}\label{EE-OE-THM-combinatorial-1}
 \sum_{\pi\in\mathcal{P}(n)}(-1)^{\ell_{N>O}(\pi)}=\sum_{\pi\in\mathcal{P}_e(n)}1-\sum_{\pi\in\mathcal{P}_o(n)}1.
\end{equation}

For an overpartition $\pi\in\mathcal{H}_{NO}(n)$, we have $SO(\pi)=LO(\pi)\geq LN(\pi)$, which yields $\ell_{N>O}(\pi)=0$, and so $(-1)^{\ell_{N> O}(\pi)}=1$.  Then, we get
\begin{equation}\label{EE-OE-THM-combinatorial-2}
 \sum_{\pi\in\mathcal{H}_{NO}(n)}(-1)^{\ell_{N>O}(\pi)}=\sum_{\pi\in\mathcal{H}_{NO}(n)}1.
\end{equation}

By restricting the involution $\mathcal{I}$ defined in Definition \ref{defi-involution} on $\mathcal{H}_{NO}(n)\bigcup\mathcal{P}(n)$, we find that the map $\mathcal{I}$ is a bijection between $\mathcal{H}_{NO}(n)$ and $\mathcal{P}(n)$, and so we have
\begin{equation}\label{EE-OE-THM-combinatorial-3}
\sum_{\pi\in\mathcal{H}_{NO}(n)}1=\sum_{\pi\in\mathcal{P}(n)}1=\sum_{\pi\in\mathcal{P}_e(n)}1+\sum_{\pi\in\mathcal{P}_o(n)}1.
\end{equation}

Using Lemma \ref{lem-involution-CF}, we have
\begin{equation*}
 \sum_{\pi\in\mathcal{C}_{NO}(n)\bigcup\mathcal{F}_{NO}(n)}(-1)^{\ell_{N> O}(\pi)}=0.
\end{equation*}
Combining with \eqref{EE-OE-THM-combinatorial-1}, \eqref{EE-OE-THM-combinatorial-2} and \eqref{EE-OE-THM-combinatorial-3}, we arrive at \eqref{EE-OE-THM-combinatorial-0}. The proof is complete. \qed

\section{Proofs of the results in Section 2.2}

In this section, we aim to prove  the results in Section 2.2. We will prove \eqref{gen-G-N-O} and \eqref{gen-E-N-O} in Section 5.1, \eqref{gen-G-N-O-new} and  \eqref{gen-E-N-O-new} in Section 5.2, Corollary \ref{equiv-g-sigma} in Section 5.3, \eqref{gen-E-N=O} and Corollary \ref{core-NO=R} in Section 5.4, Theorem \ref{N-O-THM-1} in Section 5.5, and Theorem \ref{thm-e-o-gen-1} in Section 5.6.

\subsection{Proofs of \eqref{gen-G-N-O} and \eqref{gen-E-N-O}}

In this subsection, we give the proofs of \eqref{gen-G-N-O} and \eqref{gen-E-N-O} in view of separable overpartition classes.

{\noindent \bf Proof of \eqref{gen-G-N-O}.} For $m\geq 1$, define
\[\mathcal{BG}^N_O(m)=\left\{(1^m),(\overline 1,2^{m-1}),(\overline 1,\overline 2,3^{m-2}),\ldots,(\overline 1,\overline 2,\overline 3,\ldots,\overline{m-1}, m)\right\}.\]
Set
  \[\mathcal{BG}^N_O=\bigcup_{m\geq 1}\mathcal{BG}^N_O(m).\]
Then, it is clear that $\mathcal{G}^N_O$ is a separable overpartition class and $\mathcal{BG}^N_O$ is the basis of $\mathcal{G}^N_O$. So, we get
\begin{align*}
\sum_{\pi\in\mathcal{G}_O^{N}}q^{|\pi|}&=\sum_{m\geq1}\frac{1}{(q;q)_m}\sum_{k=0}^{m-1}q^{1+2+\cdots+k+(k+1)(m-k)}\\
&=\sum_{m\geq1}\frac{1}{(q;q)_m}\sum_{k=0}^{m-1}q^{(k+1)m-\binom{k+1}{2}}.
\end{align*}
The proof is complete.  \qed

%To give a proof of \eqref{gen-E-N-O}, in the remaining of this subsection, we impose the following order on the parts of an overpartition:
%\begin{equation*}
%\overline{1}<{1}<\overline{2}<{2}<\cdots.
%\end{equation*}
%That is to say, an overpartition is a partition such that the last occurrence of a part can be overlined. For an overpartition $\pi$, we write $\pi=\left(\overline{1}^{f_{\overline{1}}(\pi)}  1^{f_1(\pi)} \overline{2}^{f_{\overline{2}}(\pi)} 2^{f_2(\pi)}\cdots\right)$. Now, we are in a position to show \eqref{gen-E-N-O}.

{\noindent \bf Proof of \eqref{gen-E-N-O}.} For $m\geq 1$, define
\[\mathcal{BE}^N_O(m)=\left\{(1^m),(1^{m-1},\overline 1),(\overline 1,2^{m-2},\overline 2),\ldots,(\overline 1,\overline 2,\ldots,m-1,\overline{m-1})\right\}.\]
Set
  \[\mathcal{BE}^N_O=\bigcup_{m\geq 1}\mathcal{BE}^N_O(m).\]
Then, it is clear that $\mathcal{E}^N_O$ is a separable overpartition class and $\mathcal{BE}^N_O$ is the basis of $\mathcal{E}^N_O$. So, we get
\begin{align*}
\sum_{\pi\in\mathcal{E}_O^{N}}q^{|\pi|}&=\sum_{m\geq1}\frac{1}{(q;q)_m}\left[q^m+\sum_{k=1}^{m-1}q^{1+2+\cdots+k+k(m-k)}\right]\\
&=\sum_{m\geq1}\frac{1}{(q;q)_m}\left[q^m+\sum_{k=1}^{m-1}q^{km-\binom{k}{2}}\right].
\end{align*}
This completes the proof.  \qed

\subsection{Proofs of \eqref{gen-G-N-O-new} and  \eqref{gen-E-N-O-new}}

In this subsection, we will show \eqref{gen-G-N-O-new} and  \eqref{gen-E-N-O-new} by considering the smallest non-overlined part.

{\noindent \bf Proof of \eqref{gen-G-N-O-new}.} In view of the smallest non-overlined part, we can get
\begin{align*}
\sum_{\pi\in\mathcal{G}_O^{N}}q^{|\pi|}&=\sum_{n\geq1}\frac{q^n}{(q^n;q)_\infty}(-q;q)_{n-1}\\
&=\frac{1}{(q;q)_\infty}\sum_{n\geq 1}q^n(q^2;q^2)_{n-1}.\\
\end{align*}
Then, it remains to show that
\begin{equation}\label{proof-GNO-1}
\frac{1}{(q;q)_\infty}\sum_{n\geq 1}q^n(q^2;q^2)_{n-1}=(-q;q)_\infty\sum_{n\geq0}\frac{q^{2n+1}}{1-q^{2n+1}}\frac{1}{(q^2;q^2)_n}.
\end{equation}

It is easy to see that
\begin{equation}\label{proof-GNO-2}
\begin{split}
\frac{1}{(q;q)_\infty}\sum_{n\geq 1}q^n(q^2;q^2)_{n-1}&=\frac{1}{(q;q)_\infty}\sum_{n\geq0}q^{n+1}(q^2;q^2)_{n}\\
&=\frac{q}{(q;q)_\infty}\sum_{n\geq0} \frac{(q^2;q^2)_n(q^2;q^2)_n}{(q^2;q^2)_n(0;q^2)_n}q^n.
\end{split}
\end{equation}

Letting $q\rightarrow q^2$, $a=b=q^2$, $c=0$ and $t=q$ in \eqref{Heine}, we get
\begin{align*}\label{proof-GNO-3}
\sum_{n\geq0} \frac{(q^2;q^2)_n(q^2;q^2)_n}{(q^2;q^2)_n(0;q^2)_n}q^n&=\frac{(q^2;q^2)_\infty(q^3;q^2)_\infty}{(0;q^2)_\infty(q;q^2)_\infty}
\sum_{n\geq0}\frac{(0;q^2)_n(q;q^2)_n}{(q^2;q^2)_n(q^3;q^2)_n}q^{2n}\\
&=(q^2;q^2)_\infty\sum_{n\geq0}\frac{q^{2n}}{(q^2;q^2)_n(1-q^{2n+1})}.
\end{align*}
Combining with \eqref{proof-GNO-2}, we have
\begin{align*}
\frac{1}{(q;q)_\infty}\sum_{n\geq 1}q^n(q^2;q^2)_{n-1}&=\frac{q}{(q;q)_\infty}(q^2;q^2)_\infty\sum_{n\geq0}\frac{q^{2n}}{(q^2;q^2)_n(1-q^{2n+1})}\\
&=(-q;q)_\infty\sum_{n\geq0}\frac{q^{2n+1}}{1-q^{2n+1}}\frac{1}{(q^2;q^2)_n}.
\end{align*}
So, \eqref{proof-GNO-1} is valid. The proof is complete. \qed

{\noindent \bf Proof of \eqref{gen-E-N-O-new}.} In virtue of the smallest non-overlined part, we can get
\begin{align*}
\sum_{\pi\in\mathcal{E}_O^{N}}q^{|\pi|}&=\sum_{n\geq1}\frac{q^n}{(q^n;q)_\infty}(-q;q)_{n}\\
&=\frac{1}{(q;q)_\infty}\sum_{n\geq 1}q^n(q^2;q^2)_{n-1}(1+q^n)\\
&=\frac{1}{(q;q)_\infty}\sum_{n\geq 1}q^n(q^2;q^2)_{n-1}+\frac{1}{(q;q)_\infty}\sum_{n\geq 1}q^{2n}(q^2;q^2)_{n-1}.
\end{align*}
Appealing to \eqref{proof-GNO-1}, we find that it remains to show
\begin{equation}\label{proof-GON-1}
\frac{1}{(q;q)_\infty}\sum_{n\geq 1}q^{2n}(q^2;q^2)_{n-1}=\frac{1}{(q;q)_\infty}-(-q;q)_\infty.
\end{equation}

Letting $q\rightarrow q^2$ and $t=-1$ in \eqref{gen-distinct}, we can get
\[\sum_{n\geq 1}q^{2n}(q^2;q^2)_{n-1}=1-(q^2;q^2)_\infty.\]

It yields that
\[\frac{1}{(q;q)_\infty}\sum_{n\geq 1}q^{2n}(q^2;q^2)_{n-1}=\frac{1}{(q;q)_\infty}\left(1-(q^2;q^2)_\infty\right)=\frac{1}{(q;q)_\infty}-(-q;q)_\infty.\]
We arrive at \eqref{proof-GON-1}. The proof is complete.  \qed

\subsection{Combinatorial proof of Corollary \ref{equiv-g-sigma}}

In this subsection, we will present the combinatorial proof of Corollary \ref{equiv-g-sigma}. For $n\geq 1$ and $k\geq 1$, we define \[\mathcal{N}(n,k)=\left\{\mu\mid|\mu|=n,0\leq k<LN(\mu)-maex(\mu)\right\}.\]
Let $\mathcal{G}_O^{N}(n,k)$ denote the set of overpartitions of $n$ in $\mathcal{G}_O^{N}$ with exactly $k$ overlined parts. It suffices to show the following lemma.

\begin{lem}\label{equiv-N-sub}
For $n\geq 1$ and $k\geq 0$, there is a bijection between $\mathcal{G}_O^{N}(n,k)$ and $\mathcal{N}(n,k)$.
\end{lem}

\pf Let $\pi=(\pi_1,\pi_2,\ldots,\pi_m)$ be a partition in $\mathcal{G}_O^{N}(n,k)$. Using the definition of  $\mathcal{G}_O^{N}(n,k)$, we deduce that $m>k$, $\pi_1,\ldots,\pi_{m-k}$ are non-overlined parts, $\pi_{m-k+1},\ldots,\pi_m$ are overlined parts and
\[|\pi_1|\geq\cdots\geq|\pi_{m-k}|>|\pi_{m-k+1}|>\cdots>|\pi_m|.\]
We first change the overlined parts $\pi_{m-k+1},\ldots,\pi_m$ in $\pi$ to non-overlined parts and denote the resulting partition by $\lambda$. Clearly, we have
\begin{equation}\label{conjugate-pro}
\lambda_1\geq\cdots\geq\lambda_{m-k}>\lambda_{m-k+1}>\cdots>\lambda_m.
\end{equation}
Then, we set $\mu$ be the conjugate of $\lambda$. We proceed to show that $\mu$ is a partition in $\mathcal{N}(n,k)$. We consider the largest part of $\mu$ and the maximal excludant of $\mu$.

It follows from the definition of conjugate that the largest part of $\mu$ is $m$, that is,
\begin{equation}\label{conjugate-largest}
LN(\mu)=m.
\end{equation}
Again by the definition of conjugate, we have  $f_t(\mu)=\lambda_{t}-\lambda_{t+1}$ for $1\leq t<m$. Appealing to \eqref{conjugate-pro}, we obtain that for $m-k\leq t<m$,
$f_t(\mu)=\lambda_{t}-\lambda_{t+1}>0$, which implies that $m-k,m-k+1,\ldots,m-1$ occur in $\mu$. Recall that
the largest part of $\mu$ is $m$, so we get
 \begin{equation}\label{conjugate-maex}
  maex(\mu)<m-k.
 \end{equation}

 Combining \eqref{conjugate-largest} and \eqref{conjugate-maex}, we have
 \[LN(\mu)-maex(\mu)>m-(m-k)=k,\]
  and so $\mu$ is a partition in $\mathcal{N}(n,k)$. Obviously,
the process above is reversible. The proof is complete. \qed

\subsection{Proofs of \eqref{gen-E-N=O} and Corollary \ref{core-NO=R}}

In this subsection, we first give an analytic proof of \eqref{gen-E-N=O} and then we give a combinatorial proof of Corollary \ref{core-NO=R}.

{\noindent \bf Analytic proof of \eqref{gen-E-N=O}.} By considering the size of the smallest non-overlined part and the largest overlined part, we can get
\begin{equation*}
\sum_{n\geq 1}{E}_{NO}(n)q^{n}=\sum_{n\geq 1}\frac{q^n}{(q^n;q)_\infty}q^{n}(-q;q)_{n-1}=\frac{1}{(q;q)_\infty}\sum_{n\geq 1}q^{2n}(q^2;q^2)_{n-1}.
\end{equation*}
Combining with \eqref{proof-GON-1}, we arrive at \eqref{gen-E-N=O}, and thus the proof is complete. \qed

Then, we give a combinatorial proof of Corollary \ref{core-NO=R}, which can be seen as a bijective proof of \eqref{gen-E-N=O}.

{\noindent \bf Combinatorial proof of Corollary \ref{core-NO=R}.} Let $\pi$ be an overpartition counted by ${E}_{NO}(n)$. We define $\lambda=\phi_{NO}(\pi)$ by changing the overlined parts in $\pi$ to non-overlined parts. Clearly, $\lambda$ is a partition enumerated by $R(n)$.

Conversely, let $\lambda$ be a partition enumerated by $R(n)$. Assume that $k$ is the smallest integer such that $k$ appears at least twice in $\lambda$. Then, we define $\pi=\psi_{NO}(\lambda)$ by changing a $k$ in $\lambda$ to $\overline k$ and changing the parts smaller than $k$ in $\lambda$ to overlined parts. It is obvious that $\pi$ is an overpartition counted by ${E}_{NO}(n)$. Clearly, $\psi_{NO}$ is the inverse map of $\phi_{NO}$. Thus, we complete the proof.  \qed

\subsection{Proofs of Theorem \ref{N-O-THM-1}}

This subsection is devoted to presenting an analytic proof and a combinatorial of Theorem  \ref{N-O-THM-1}. We first give the analytic proof of Theorem  \ref{N-O-THM-1}.

{\noindent \bf Analytic proof of Theorem \ref{N-O-THM-1}.} In virtue of the smallest non-overlined part, we can get
\begin{equation*}
\sum_{n\geq 1}\left({A}_{O\geq N}(n)-{B}_{O\geq N}(n)\right)q^n=\sum_{n\geq 1} \frac{q^n}{(q^n;q)_\infty}(-q;q)_{n-1}(q^n;q)_\infty
=\sum_{n\geq1}q^n(-q;q)_{n-1}.
\end{equation*}

Setting $t=1$ in \eqref{gen-distinct}, we find that
\[\sum_{n\geq1}q^n(-q;q)_{n-1}\]
is the generating function for the nonempty distinct partitions.
The proof is complete.  \qed

Before giving the combinatorial proof of Theorem \ref{N-O-THM-1}, for $n\geq 1$, we introduce the following notations.
\begin{itemize}

\item Let $\overline{\mathcal{PN}}(n)$ be the set of overpartitions of $n$ such that non-overlined parts must appear.

\item Let $\mathcal{C}_{ON}(n)$ be the set of overpartitions $\pi$ of $n$ with $LO(\pi)\geq LN(\pi)\geq1$.

\item Let $\mathcal{F}_{ON}(n)$ be the set of overpartitions $\pi$ of $n$ such that there are at least two non-overlined parts in $\pi$ and $LN(\pi)> LO(\pi)$.

\item Let $\mathcal{H}_{ON}(n)$ be the set of overpartitions $\pi$ of $n$ such that there is exactly one non-overlined part in $\pi$ and $LN(\pi)> LO(\pi)$.

\end{itemize}
Clearly, we have
\[\overline{\mathcal{PN}}(n)=\mathcal{C}_{ON}(n)\bigcup\mathcal{F}_{ON}(n)\bigcup\mathcal{H}_{ON}(n).\]

By restricting the involution $\mathcal{I}$ defined in Definition \ref{defi-involution} on $\mathcal{C}_{ON}(n)\bigcup\mathcal{F}_{ON}(n)$, it is easy to get the following lemma.
\begin{lem}\label{lem-involution-CF-NEXT}
For $n\geq 1$, the map $\mathcal{I}$ is a bijection between $\mathcal{C}_{ON}(n)$ and $\mathcal{F}_{ON}(n)$. Moreover,  for an overpartition $\pi \in \mathcal{C}_{ON}(n)$, let $\lambda=\mathcal{I}(\pi)$. We have $\ell_{O\geq N}(\lambda)=\ell_{O\geq N}(\pi)-1$.
\end{lem}

With a bijective method, we get
\begin{lem}\label{lem-ON-H-NUMBER-11111}
For $n\geq 1$, the number of overpartitions in $\mathcal{H}_{ON}(n)$ is $D(n)$.
\end{lem}

\pf Let $\pi=(\pi_1,\pi_2,\ldots,\pi_m)$ be an overpartition in $\mathcal{H}_{ON}(n)$. If we change the overlined parts $\pi_2,\ldots,\pi_m$ in $\pi$ to non-overlined parts, then we get a distinct partition counted by $D(n)$, and vice versa. This completes the proof. \qed

Now, we are in a position to give the combinatorial proof of Theorem \ref{N-O-THM-1}.

{\noindent \bf Combinatorial proof of Theorem \ref{N-O-THM-1}.}
Clearly, it is equivalent to showing that for $n\geq 1$,
 \begin{equation}\label{EO-OO-THM-combinatorial-0}
 \sum_{\pi\in\overline{\mathcal{PN}}(n)}(-1)^{\ell_{O\geq N}(\pi)}=D(n).
\end{equation}

Appealing to Lemma \ref{lem-involution-CF-NEXT}, we have
\begin{equation}\label{EO-OO-THM-combinatorial-1}
 \sum_{\pi\in\mathcal{C}_{ON}(n)\bigcup\mathcal{F}_{ON}(n)}(-1)^{\ell_{O\geq N}(\pi)}=0.
\end{equation}

For an overpartition $\pi\in\mathcal{H}_{ON}(n)$, we have $SN(\pi)=LN(\pi)> LO(\pi)$, which yields $\ell_{O\geq N}(\pi)=0$, and so $(-1)^{\ell_{O\geq N}(\pi)}=1$. Using Lemma \ref{lem-ON-H-NUMBER-11111}, we get
\begin{equation*}
 \sum_{\pi\in\mathcal{H}_{ON}(n)}(-1)^{\ell_{O\geq N}(\pi)}=D(n).
\end{equation*}
Combining with \eqref{EO-OO-THM-combinatorial-1}, we arrive at \eqref{EO-OO-THM-combinatorial-0}. The proof is complete.  \qed

\subsection{Proofs of Theorem \ref{thm-e-o-gen-1}}

In this subsection, we will give an analytic proof and a combinatorial of Theorem \ref{thm-e-o-gen-1}. We first present the analytic proof of Theorem \ref{thm-e-o-gen-1}.

{\noindent \bf Analytic proof of Theorem \ref{thm-e-o-gen-1}.} It is clear that
\begin{equation}\label{final-000000000}
\sum_{n\geq 1}{H}'_{ON}(n)q^n=\sum_{n\geq1}
\frac{q^n}{1-q^n}(-q;q)_n.
\end{equation}

In view of the smallest non-overlined part, we can get
\begin{align*}
\sum_{n\geq 1}\left({A}_{O> N}(n)-{B}_{O> N}(n)\right)q^n&=\sum_{n\geq1}\frac{q^n}{(q^n;q)_\infty}(-q;q)_n(q^{n+1};q)_\infty\\
&=\sum_{n\geq1}\frac{q^n}{1-q^n}(-q;q)_n.
\end{align*}
Combining with \eqref{final-000000000}, we complete the proof.  \qed

Then, we give the combinatorial of Theorem \ref{thm-e-o-gen-1} with a similar argument in the combinatorial proof of Theorem \ref{N-O-THM-1}.

{\noindent \bf Combinatorial proof of Theorem \ref{thm-e-o-gen-1}.} Clearly, it is equivalent to showing that for $n\geq 1$,
\begin{equation}\label{EO-OO-THM-combinatorial-000-0}
 \sum_{\pi\in\overline{\mathcal{PN}}(n)}(-1)^{\ell_{O>N}(\pi)}={H}'_{ON}(n).
\end{equation}

To do this, we introduce the following three sets.
\begin{itemize}
\item Let $\mathcal{C}'_{ON}(n)$ be the set of overpartitions $\pi$ of $n$ with $LO(\pi)\geq LN(\pi)\geq 1$ and $LO(\pi)>SN(\pi)$.

\item Let $\mathcal{F}'_{ON}(n)$ be the set of overpartitions $\pi$ of $n$ with $LN(\pi)> LO(\pi)$ and $LN(\pi)> SN(\pi)\geq 1$.

\item Let $\mathcal{H}'_{ON}(n)$ be the set of overpartitions $\pi$ of $n$ with $LN(\pi)=SN(\pi)\geq 1$ and $SN(\pi)\geq LO(\pi)$.

\end{itemize}
Clearly, we have
\[\overline{\mathcal{PN}}(n)=\mathcal{C}'_{ON}(n)\bigcup\mathcal{F}'_{ON}(n)\bigcup\mathcal{H}'_{ON}(n).\]

By definition, we know that the number of overpartitions in $\mathcal{H}'_{ON}(n)$ is ${H}'_{ON}(n)$.
For an overpartition $\pi\in\mathcal{H}'_{ON}(n)$, we have $SN(\pi)\geq LO(\pi)$, which yields $\ell_{O> N}(\pi)=0$, and so $(-1)^{\ell_{O> N}(\pi)}=1$. So, we get
\begin{equation}\label{EO-OO-THM-combinatorial-000-1}
 \sum_{\pi\in\mathcal{H}'_{ON}(n)}(-1)^{\ell_{O> N}(\pi)}={H}'_{ON}(n).
\end{equation}

By restricting the involution $\mathcal{I}$ defined in Definition \ref{defi-involution} on $\mathcal{C}'_{ON}(n)\bigcup\mathcal{F}'_{ON}(n)$, we find that the map $\mathcal{I}$ is a bijection between $\mathcal{C}'_{ON}(n)$ and $\mathcal{F}'_{ON}(n)$. Moreover,  for an overpartition $\pi \in \mathcal{C}'_{ON}(n)$, let $\lambda=\mathcal{I}(\pi)$. We have $\ell_{O> N}(\lambda)=\ell_{O> N}(\pi)-1$. This implies that
\begin{equation*}
 \sum_{\pi\in\mathcal{C}'_{ON}(n)\bigcup\mathcal{F}'_{ON}(n)}(-1)^{\ell_{O> N}(\pi)}=0.
\end{equation*}
Combining with \eqref{EO-OO-THM-combinatorial-000-1}, we arrive at \eqref{EO-OO-THM-combinatorial-000-0}, and thus the proof is complete.  \qed

\noindent{\bf Acknowledgments.} This work
was supported by   the National Natural Science Foundation of China and Sichuan Science and Technology Program (No. 2024NSFSC1394).

\noindent{\bf Declarations}

\noindent{\bf Conflict of interest:} The authors have no competing interests to declare that are
relevant to the content of this article.

\end{document}